\documentclass[12pt]{article}
\usepackage{amsmath}
\usepackage{amsthm}
\usepackage{amsfonts}

\usepackage{amssymb}
\usepackage{graphicx}%
\usepackage[latin1]{inputenc}
\providecommand{\U}[1]{\protect\rule{.1in}{.1in}}

\protect \hoffset -1 true cm \voffset -1.5 true cm \textwidth 15cm
\newcommand{\ieq}{\begin{equation}}
\newcommand{\eeq}{\end{equation}}
\newcommand{\ieqa}{\begin{eqnarray}}
\newcommand{\eeqa}{\end{eqnarray}}
\newcommand{\ieqas}{\begin{eqnarray*}}
\newcommand{\eeqas}{\end{eqnarray*}}

\newtheorem{theorem}{Theorem}

\newtheorem{prop}{Proposition}

\begin{document}

\author{A. Alvino $^{(\ast)}$
\and R. Volpicelli ${^{(\ast)}}$
\and B. Volzone ${^{(\ast\ast)}}$
\and $^{(\ast)}$ Universit\`{a} degli Studi di Napoli ``Federico II"
\and Dipartimento di Matematica e Applicazioni ``Renato Caccioppoli"
\and emails: angelo.alvino@unina.it; rvolpice@unina.it
\and $^{(\ast\ast)}$ Universit\`{a} degli Studi di Napoli ``Parthenope"
\and Dipartimento per le Tecnologie
\and email: bruno.volzone@uniparthenope.it}
\title{On Hardy inequalities with a remainder term
}
\date{}
\maketitle

\begin{abstract}
In this paper we study some improvements of the classical Hardy inequality. We add to the right hand side of the inequality a term which depends on some Lorentz norms of $u$ or of its gradient and we  find the best values of the constants for remaining terms. In both cases we show that the problem of finding the optimal value of the constant can be reduced to a spherically symmetric situation. This result is new when the right hand side is a Lorentz norm of the gradient.

\end{abstract}

\date{}
\maketitle

\bigskip
{\normalsize \parindent0pt }

\section{{\protect\normalsize Introduction}}

The classical Hardy inequality asserts that (see \cite{Hardy1} and \cite{HLP})%
\begin{equation}
\int_{\Omega}\left\vert \nabla u\right\vert ^{2}dx\geq\frac{\left(
N-2\right)  ^{2}}{4}\int_{\Omega}\frac{u^{2}}{\left\vert x\right\vert ^{2}%
}dx\text{, \ \ }\forall u\in H_{0}^{1}\left(  \Omega\right)  , \label{eq.4}%
\end{equation}
where $\Omega$ is a bounded open set of $\mathbb{R}^{N}$
containing the origin, $N>2$. The constant in (\ref {eq.4}) is the best possible; however  it is not attained. This fact allows to add to the right hand side of (\ref{eq.4})  a suitable remaining term involving some
norm of $u$ or of the gradient of $u$.
The first result in this direction was obtained by Brezis and
Vazquez in \cite{Brezis 2};  if $N>2$ they show  the so called
\textit{Hardy-Poincar\`{e} inequality}%
\begin{equation}
\int_{\Omega}\left\vert \nabla u\right\vert ^{2}dx-\frac{\left(
N-2\right) }{4}^{2}\int_{\Omega}\frac{u^{2}}{\left\vert
x\right\vert ^{2}}dx\geq
\frac{ \Lambda_{2}}{R_\Omega^2}\left\Vert u\right\Vert _{2}^{2}\text{, \ }\forall u\in
H_{0}^{1}\left(
\Omega\right)  \label{eq.5}%
\end{equation}
where $\Lambda_{2}$ denotes the first eigenvalue of the Laplace
operator in the  unit disk  of $\mathbb{R}^{2}$ and $R_\Omega$ is the radius of the ball $\Omega^\#\subseteq \mathbb{R}^{N}$ centered at the origin having the same measure as $\Omega$.  The constant in (\ref  {eq.5})
 is  \textit{optimal} even if it is again not achieved.

The first aim of this paper is to find the best value of the constant $C$ in inequalities of the type (\ref {eq.5}) that involve a Lorentz norm of $u$ as a remainder term. In particular we focus our attention to the following two inequalities:
\begin{equation}
\int_{\Omega}|\nabla u|^{2}dx-\frac{(N-2)^{2}}{4}\int_{\Omega}\frac{u^{2}%
}{|x|^{2}}dx\geq C(|\Omega|) \left\Vert u\right\Vert^2_{\frac{2N}{N-1},2}\text{, \ }\forall u\in
H_{0}^{1}\left(
\Omega\right)
\label{tipogenerale1}%
\end{equation}%
\begin{equation}
\int_{\Omega}|\nabla u|^{2}dx-\frac{(N-2)^{2}}{4}\int_{\Omega}\frac{u^{2}%
}{|x|^{2}}\,dx\geq C(p,|\Omega|)\left\Vert u\right\Vert^2_{p,1}\text{, \ }\forall u\in
H_{0}^{1}\left(
\Omega\right) \label{caso31}%
\end{equation}
where $N>2$ and $1 \leq p<2^{*}$, $2^{*}=2N/(N-2)$. We recall that $u$ belongs to the Lorentz space $L\left(r,s\right)\left(\Omega\right)$ with $0<r,s<+\infty$, if the quantity defined by
\[
\left\Vert u\right\Vert_{r,s}=\omega_{N}^{1/r-1/s}\left(\int_{\Omega^{\#}} \left[u^{\#}\left(x\right)\left\vert x\right\vert^{N/r}\right]^{s}\frac{dx}{\left\vert x\right\vert^{N}}\right )^{1/s}
\]
is finite, where $u^{\#}$ is the \textit{spherical decreasing rearrangement} of $u$, which is a spherically symmetric function defined on $\Omega^\#$,  decreasing along the radius, having the same distribution function as $u$.
Inequalities (\ref{tipogenerale1}) and (\ref{caso31})  are  not new (see \cite {Chauduri}, \cite{Ghoussoub 1}). At least for (\ref{tipogenerale1}), the best value of the constant might be also obtained using the results of \cite{Ghoussoub 1}.
The technique used to get (\ref{eq.5}) (or  any of the inequalities we quoted before) follows a usual procedure. First we observe that we can restrict our attention to spherically symmetric functions defined on $\Omega^\#$.  This can be done replacing $u$ by $u^{\#}$. Indeed,  this operation decreases the left hand side of inequality (\ref{eq.5}), since it  it decreases  the $H_0^1(\Omega)$-norm  of $u$ by the classical  Pólya-Szegö  principle (see \cite{Polya}) and, by Hardy-Littlewood inequality (see \cite {HLP}, \cite {Bennet}), it increases the weighted $L^{2}$-norm of $u$ with weight $\left\vert x\right\vert^{-2}$.
Moreover, this operation does not change the
norm of $u$ in  Lebesgue spaces or, more in general  in  Lorentz spaces . Once reduced the problem to spherically symmetric case, the best value of the constant in (\ref{eq.5}) is obtained by using what  Brezis and Vazquez (see \cite{Brezis 2}) call the \textit{magical transformation}
\begin{equation}
v(r)=u(r) r^{\frac{N-2}{2}}\qquad r=|x|. \label{sostituzione}%
\end{equation}
This transformation produces a  dimension reduction  of the problem from $N$ to $2$ dimensions which also explain the presence of the constant $\Lambda_2$ in the inequality (\ref{eq.5}).

In view of inequality (\ref{eq.5}), one can also wonder if it may
be possible to replace the $L^{2}$-norm of $u$ by a norm of the gradient $\nabla
u$ of $u$. This question was showed to have a positive
answer too. Indeed, Vazquez and Zuazua (see
\cite{Vazquez1} ) proved the following \textit{Improved Hardy-Poincar\`{e} inequality}%

\begin{equation}
\int_{\Omega}\left\vert \nabla u\right\vert ^{2}dx-\frac{\left(
N-2\right) }{4}^{2}\int_{\Omega}\frac{u^{2}}{\left\vert
x\right\vert ^{2}}dx\geq
C(q,|\Omega|)\left\Vert \left\vert\nabla u\right\vert\right \Vert _{q} ^{2}\text{, \ \ }\forall u\in H_{0}%
^{1}\left(  \Omega\right) ,  \label{eq.20}%
\end{equation}
where $\Omega$ is a bounded open set of $\mathbb{R}^{N}$ containing the origin, $N>2$, $1\leq q<2$ and $C(q,\left \vert \Omega \right \vert)$ is a constant depending only on
$q$ and $\Omega$. However the question about finding the best
value of the constant $C(q,\left\vert\Omega\right\vert)$ seems to be still open and it
is not even known if it is possible, in order to find this best
value, to restrict the attention to spherically symmetric
functions defined on a ball. As regards this last question we  prove that as for inequality (\ref{eq.20}) it is possible to reduce the problem to a spherically symmetric case using a suitable symmetrization procedure. In this situation, unlike the cases explained before, it is not useful to replace the function $u$ by  its spherical decreasing rearrangement
$u^{\#}$ since  by Hardy-Littlewood inequality and by
Pólya-Szegö principle   both terms of the
inequality (\ref{eq.20})  decrease under spherical symmetrization. Our idea
is then to fix not the rearrangement of $u$ but the rearrangement
of $\left\vert\nabla u\right\vert$. In this way both
the $L^{2}$-norm and the $L^q$-norm of $\left\vert\nabla
u\right\vert$ do not change and we have only to investigate what
happens to the $L^{2}$-norm of $u$ with weight $\left\vert
x\right\vert^{-2}$. It can be seen that if $u\in H_0^1(\Omega)$ there exists a spherically symmetric function
$\overline{u}$ defined on the ball $\Omega^{\#}$, such that $\left\vert\nabla
\overline{u}\right\vert^{\#}=\left\vert\nabla
{u}\right\vert^{\#}$ and
$\int_{\Omega}\frac{u^{2}}{\left\vert x\right\vert^2}dx$ \textit{increases} when we pass from $u$ to $\overline{u}$. That will be enough for our purpose. Afterwards, we find the best value of the constant $C$ in the inequality (\ref{eq.20})  in
the case $q=1$. Indeed in this case it is quite easy to see that, if $u$ is a
spherically symmetric function defined on a ball, the problem of finding the best value of the constant in the inequality  (\ref{eq.20}) can be reduced to the study of
the  inequality  (\ref{caso31}).\\It is  clear that the same arguments apply when the
$L^{q}$-norm of $\nabla u$ is replaced by a more general Lorentz
norm $\left\Vert\left\vert\nabla u\right\vert\right\Vert_{p,q}$
with $1\leq q\leq p<2$ (see section 3).

\section{{\protect\normalsize Best constant in Hardy- Sobolev inequalities
with a remainder term in Lorentz spaces}}

In this section we prove inequalities (\ref{tipogenerale1}) and (\ref{caso31}). More precisely, we have the following

\begin{theorem}
Let $\Omega$ be a bounded open set of
$\mathbb{\mathbb{R}}^{N}$ containing the origin, $N>2$.  The optimal value
of the constant in the  inequality \textrm{(\ref{tipogenerale1})} is given by
\begin{equation}
C(|\Omega|)= \frac{\omega^{\frac {2}{N}}_{N}}{|\Omega|^{\frac{1}{N}}}V_{0}, \label{eq.7}
\end{equation}
where ${\omega_{N}}$ is the measure of the N-dimensional unit ball
and $V_{0}$ is the first zero of the function
$V(r)=J_{0}\left(2\sqrt{r}\right)$  (here $J_{0}$ denotes, as usual, the Bessel function of zero order).%

\end{theorem}

{\normalsize \noindent\textbf{Proof .} As pointed out in the introduction it is enough to prove inequality (\ref{tipogenerale1}) for spherically symmetric and decreasing functions defined on the ball $\Omega^{\#}$. }{\normalsize Hence we are reduced to study the inequality
\begin{equation}
\int_{\Omega^{\#}}|\nabla u|^{2}dx-\frac{(N-2)^{2}}{4} \int_{\Omega^{\#}}\frac{u^{2}%
}{|x|^{2}}dx\geq C\int_{\Omega^{\#}}\frac{u^{2}}{|x|}dx\qquad u\in
H_{0}^{1}({\Omega^{\#}})
\label{tipoparticolare3}%
\end{equation}
where  $u$ is a radial function.  For simplicity we will study the inequality
\begin{equation}
\int_{B_R}|\nabla u|^{2}dx-\frac{(N-2)^{2}}{4} \int_{B_R}\frac{u^{2}%
}{|x|^{2}}dx\geq \int_{B_R}\frac{u^{2}}{|x|}dx\qquad u\in
H_{0}^{1}(B_R)
\label{tipoparticolare2}%
\end{equation}
where $B_R$ is the ball centered at the origin whose radius $R$ must be determined and $u$ is a radial function. This is not a restriction since  (\ref {tipoparticolare3}) reduces to (\ref {tipoparticolare2}) using a suitable homothety. Indeed if $u=u(|x|)$ is a radial function for which  (\ref {tipoparticolare3}) holds the function $z=u\left(\frac {|x|} {C}\right)$ satisfies (\ref {tipoparticolare2}) with $R=C\left(\frac {|\Omega|} {\omega_N}\right)^{1/N}$.
We make the classical
change of variable (\ref{sostituzione}) which essentially allows us to read the Hardy-Sobolev inequality
\textrm{(\ref{tipoparticolare2})} as a Sobolev inequality in the
plane. Indeed assuming
$u\in C_{0}^{1}(B_{R})$ and hence $v(0)=0$, inequality
(\ref{tipoparticolare2}) becomes
\begin{equation}
\int_{0}^{R}(v^{\prime})^{2}r\,dr\geq\int_{0}^{R}v^{2}\,dr\,\qquad
v\in
H^{1}\left(0,R\right),\>v(R)=0, \label{Sobolev}%
\end{equation}
Let us consider the functional
\begin{equation}
J(v)=\int_{0}^{R}(v^{\prime})^{2}r\,dr-\int_{0}^{R}v^{2}\,dr.
\end{equation}
The Euler equation of this functional is
\[
(rv^{\prime})^{\prime}+v=0\;
\]
a solution of this equation is the function
\[
V(r)=J_{0}\left(2\sqrt{r}\right)=\sum_{n=0}^{\infty}(-1)^{n}\frac{r^{n}}{(n!)^{2}}\,.
\]
By standard arguments of Calculus of Variations it can be seen
that the functions $v(r)=cV(r)$ with $c\in\mathbb{R}$ minimize the
functional $J\left(v\right)$.
Hence (\ref{Sobolev}) holds for $R=V_{0}$ where  $V_{0}\simeq1.4457\ldots$ is the first zero of the function
$V(r)=J_{0}\left(2\sqrt{r}\right)$ .

\noindent Using again the change of variable (\ref{sostituzione})
and coming back to the function $u$ we obtain
(\ref{tipoparticolare2}) for $R=V_{0}$.
The restriction $u\in C_{0}^{1}(B_{R})$
can be removed by density. A dimensional analysis on the constant
shows that   inequality (\ref{tipoparticolare3}) holds for
$C=\left(\frac{\omega_{N}}{|\Omega|}\right)^{\frac{1}{N}}V_{0}$. Hence the best value of the constant in the inequality (\ref{tipogenerale1}) is given by (\ref {eq.7}).
On the other hand it is clear that the
optimal value of the constant is not attained since it would
correspond to equality in (\ref{Sobolev}) which happens for
$v(r)=cV(r)$ and hence
\[
u(x)=c|x|^{-\frac{N-2}{2}}V(|x|)
\]
which is not in $H^{1}$.
$\qedsymbol$

An improvement of inequality  (\ref{tipogenerale1}) can be obtained by introducing  a weighted norm of $u$ with a logarithmic weight. The next proposition gives a very simple proof of this type of inequality.
These inequalities are widely studied in \cite{Filippas 2}, \cite {Filippas 3}, \cite {Filippas}.

\begin{prop}
Let $\Omega$ be a bounded open set of $\mathbb{R}^{N}$ containing the origin, $N>2$, of measure $\displaystyle\frac{\omega_{N}%
}{e^{N}}$ and let $p$  be a function defined
in $\Omega
$ whose spherical decreasing rearrangement is $p^{\#}%
(x)=\displaystyle\dfrac{1}{[|x|\log|x|]^{2}}\,.$ Then the following inequality holds
\begin{equation}
\int_{\Omega}|\nabla u|^{2}dx-\frac{(N-2)^{2}}{4}\int_{\Omega}\frac{u^{2}%
}{|x|^{2}}dx\geq\frac{1}{4}\int_{\Omega}p(x)u^{2}dx. \label{greca}%
\end{equation}
\end{prop}

{\normalsize \noindent\textbf{Proof .}
The inequality (\ref{greca})
was first proved in \cite{Brezis 1} in the case }$N=1$
and more recently in \cite{Filippas 2} in any dimension.
{\normalsize Pólya-Szegö principle and
Hardy-Littlewood inequality allows us to reduce the study of
(\ref{greca}) to radial function
defined in the ball of radius $\frac{1}{e}$. In this case (\ref{greca}) becomes %
\begin{equation}
\int_{B_{1/e}}|\nabla u|^{2}dx-\frac{(N-2)^{2}}{4}\int_{B_{1/e}}\frac{u^{2}%
}{|x|^{2}}dx\geq\frac{1}{4}\int_{B_{1/e}}\frac{u^{2}}{[|x|\log|x|]^{2}}dx
\label{grecasimm}%
\end{equation}
We just observe that making the change of variable  (\ref{sostituzione})%
$\>\>$  with $v$ such that $v(0)=v(1/e)=0$ we get
just the inequality proved in \cite{Brezis 1}. We give here an
alternative and even simpler proof of (\ref{grecasimm}). } }{\normalsize Making the change of variable %
\[
u(r)=v(r)r^{-\frac{N-2}{2}}\sqrt{-\log r}%
\]
with $v$  such that $v(1/e)=0$ and
\begin{equation}
\lim_{r\rightarrow0}v(r)\log r=0\,, \label{lim0}%
\end{equation}
it results%
\begin{align*}
\int_{0}^{1/e}(u^{\prime})^{2}r^{N-1}dr  &  =\displaystyle\int_{0}%
^{1/e}(v^{\prime})^{2}|\log r|rdr+\displaystyle\frac{\left
(N-2\right)^{2}}{4}
\displaystyle\int_{0}^{1/e}\frac{v^{2}|\log r|}{r}dr\\
& \\
&  +\frac{1}{4}\int_{0}^{1/e}\frac{v^{2}}{r|\log r|}dr-(N-2)\int_{0}%
^{1/e}vv^{\prime}|\log r| dr\\
& \\
&  -\int_{0}^{1/e}vv^{\prime}dr+\frac{N-2}{2}\int_{0}^{1/e}\frac{v^{2}}%
{r}dr\,.
\end{align*}
}

{\normalsize \noindent Using the boundary condition on
}$v$ we have %
\[
\int_{0}^{1/e}vv^{\prime}\log rdr=-\frac{1}{2}\int_{0}^{1/e}%
\frac{v^{2}}{r}dr\,.
\]
Hence%
\[
\int_{B_{1/e}}|\nabla u|^{2}dx\geq N\omega_{N}\displaystyle\frac
{\left(N-2\right)^{2}}{4}\displaystyle\int_{0}^{1/e}\frac{v^{2}|\log
r|}{r}dr+\frac {1}{4}N\omega_{N}\int_{0}^{1/e}\frac{v^{2}}{r|\log
r|}dr
\]
that is (\ref{grecasimm}).$\qedsymbol$

\bigskip

Much more interesting seems to us inequality (\ref{caso31}) since it will be used in section 3 to find the best value of the constant in inequalities involving as a remainder term the $L(p,1)$ norm of the gradient for $0<p<\infty$. We have

\begin{theorem}
Let $\Omega$ be a bounded open set of $\mathbb{R}%
$$^{N}$ containing the origin, $N>2$. Then for any 1 $\leq p<2^{*} $ , the optimal
value of the constant in the inequality (\ref {caso31}) is given by
\begin{equation}
C(p,|\Omega|)=\frac{2\left(  \frac{N}{p}-\frac{N}{2}+1\right)  ^{3}%
}{N\left\vert \Omega\right\vert^{\frac{2\left(  \frac{N}{p}-\frac{N}%
{2}+1\right)  }{N}}}\cdot\omega_N^{
\frac{2}{N}
}.\label{migcost2}
\end{equation}

\end{theorem}

{\normalsize \noindent\textbf{Proof .} As in the proof of theorem
1, using symmetrization we are reduced to prove inequality (\ref{caso31}) for radial function $u$ defined on $\Omega^\#$. Moreover if we replace $u$ by the function  $z=u\left(\frac {|x|} {K}\right)$ with $K^{2(\frac {N}{p}-\frac{N}{2}+1)}=C$ we see that it is enough to study
\begin{equation}
\int_{B_{R}}|\nabla u|^{2}dx-\frac{(N-2)^{2}}{4}\int_{B_{R}}\frac{u^{2}%
}{|x|^{2}}\,dx\geq \left(  \int_{B_{R}}\frac{|u|}{|x|^{N-\frac{N}{p}}%
}\,dx\right)  ^{2}\,\qquad u\in H_{0}^{1}(B_{R}) \label{casosimmetrico31}%
\end{equation}
where $u$ is a radial function defined in a suitable ball $B_{R}$
centered at the origin. }

{\normalsize \noindent Making the change of variable
\textrm{(\ref{sostituzione})} and assuming $u\in
C_{0}^{1}(B_{R})$, we obtain the following inequality in the plane
\begin{equation}
\int_{C_{R}}|\nabla v|^{2}dx\geq\frac{N\omega_{N}}{2\pi}\left(  \int_{C_{R}%
}\frac{v}{|x|^{1-\frac{N}{p}+\frac{N}{2}}}dx\right)  ^{2}\,
\label{stimafinale}%
\end{equation}
where $C_{R}$ is the disk centered at the origin. Let us consider
the functional
\[
I(v)=\int_{C_{R}}|\nabla
v|^{2}dx-\frac{N\omega_{N}}{\pi}\int_{C_{R}}\frac
{v}{|x|^{1-\frac{N}{p}+\frac{N}{2}}}dx=2\pi\int_{0}^{R}(v^{\prime}%
)^{2}r\,dr-\,2N\omega_{N}\int_{0}^{R}\frac{v}{r^{\frac{N}{2}-\frac{N}{p}}%
}\,dr.\
\]
The function
\[
V(r)=\frac{N\omega_{N}}{2\pi}\cdot\frac{1}{\left(  \frac{N}{p}-\frac{N}%
{2}+1\right)  ^{2}}\left(  R^{\frac{N}{p}-\frac{N}{2}+1}-r^{\frac{N}{p}%
-\frac{N}{2}+1}\right)
\]
minimizes this functional. Moreover if
\begin{equation}
R^{2\left(  \frac{N}{p}-\frac{N}{2}+1\right)  }=\frac{2\left(  \frac{N}%
{p}-\frac{N}{2}+1\right)  ^{3}}{N\omega_{N}}, \label{raggio}%
\end{equation}
it satisfies
\begin{equation}
\int_{C_{R}}\frac{V}{|x|^{1-\frac{N}{p}+\frac{N}{2}}}dx=1\,
\label{normalizzazione}%
\end{equation}
and
\begin{equation}
\int_{C_{R}}|\nabla V|^{2}dx=\frac{N\omega_{N}}{2\pi}. \label{eq.3}%
\end{equation}
}

{\normalsize \noindent If $v$ is any spherically symmetric
function in $H_{0}^{1}(C_{R})$ for which condition
\textrm{(\ref{normalizzazione})} holds, by (\ref{eq.3}) we have
\[
\int_{C_{R}}|\nabla v|^{2}dx-\frac{N\omega_{N}}{\pi}=I(v)\geq
I(V)=-\frac
{N\omega_{N}}{2\pi}%
\]
and then

\[
\int_{C_{R}}|\nabla v|^{2}dx\geq \frac
{N\omega_{N}}{2\pi}.%
\]

\noindent Hence we find
\begin{equation}
\int_{C_{R}}|\nabla v|^{2}\geq\frac{N\omega_{N}}{2\pi}\left(  \int_{C_{R}%
}\frac{v}{|x|^{1-\frac{N}{p}+\frac{N}{2}}}\right)  ^{2} \label{best}%
\end{equation}
if we remove the assumption \textrm{(\ref{normalizzazione})} on
$v$. Coming back to the function $u,$ using (\ref{best}) we get
(\ref{casosimmetrico31}) with $R$
  given in (\ref{raggio}) . The assumption $u\in C_{0}^{1}(B_{R})$
can be removed by density. A dimensional analysis on the constant
shows that inequality (\ref{caso31}) holds for
$
C(p,|\Omega|)$ given by (\ref {migcost2}).$\qedsymbol$

\section{{\protect\normalsize Best constant in Hardy-Sobolev inequalities with
a remainder term depending on the gradient}}

{\normalsize In this section we focus our attention on the
inequality (\ref{eq.20}).
As pointed out in the introduction, our aim is to show a new approach that allows us to reduce the
problem of finding the optimal value of the constant  in  (\ref{eq.20}) to a spherically
symmetric situation.  Our idea is  to fix not the rearrangement of
$u$, as we did in the proofs of the inequalities of the previous section, but the rearrangement of its gradient. In order to explain
into details the result let us recall some preliminaries. Let $f$, $g$ be two non negative measurable functions defined on two open bounded sets  having the same measure $V$.
We say that $f$ is
\textit{dominated} by $g$, and we write $f\prec g$, if%
\[
\int_{0}^{s}f^{\ast}(\sigma)d\sigma\leq\int_{0}^{s}g^{\ast}\left(
\sigma\right)  d\sigma\text{, \ }\forall s\in\left[  0,V\right]
\]
and%
\[
\int_{0}^{V}f^{\ast}(\sigma)d\sigma=\int_{0}^{V}g^{\ast}\left(
\sigma\right) d\sigma.
\]
We say that $f$ \textit{is equimeasurable with} $g,$ or that $f$
\textit{is a
rearrangement of} $g$, if%
$\>\>f^{\ast}=g^{\ast}$.

\noindent Let $f_{0}$ be a prescribed non negative decreasing and right-continuous function from $L^{p}\left(  0,V\right)$, $p\geq 1$.

\noindent It is known that (see \cite{AT}, \cite{Mossino} and also \cite{EdyRoberta} for further details) if $f\in
L^{p}\left(  0,V\right)  $ and $f\prec f_{0}$, there
exists a sequence $\left\{ f_{n}\right\}
_{n\in\mathbb{N}}\subseteq L^{p}\left(0,V\right)$ such that
$f_{n}^{\ast}=f_{0}$ for all $n\in\mathbb{N}$
and $f_{n}$ converges weakly to $f$ in  $L^{p}\left(  0,V\right)$.

\bigskip
The following theorem shows that the value of the functional
\[
J\left(  u\right)  =\frac{
{\displaystyle\int_{\Omega}}
\left\vert \nabla u\right\vert ^{2}dx-\dfrac{\left(  N-2\right)  ^{2}}{4}%
{\displaystyle\int_{\Omega}}
\dfrac{u^{2}}{\left\vert x\right\vert
^{2}}dx}{\left\Vert\left\vert\nabla u\right\vert\right\Vert
_{q}^{2}}%
\]
decreases if we replace $\Omega$ with the ball
$\Omega^{\#}$, the function $u$ by a suitable spherically
symmetric and decreasing function $\overline{u}\in H_{0}^{1}\left(
\Omega^{\#}\right)  $ such that $\left\vert
\nabla\overline{u}\right\vert $ is a rearrangement of $\left\vert
\nabla u\right\vert $.

\begin{theorem}
Let $\Omega$ be a bounded,  open set of $\mathbb{R}^{N}$ containing the origin, $N>2$.
If $u\in H_{0}^{1}\left(  \Omega\right)  $ is a non negative
function,  then there exists a  spherically symmetric decreasing function
$\overline{u}\in H_{0}^{1}\left( \Omega^{\#}\right)  $ such that
$\left\vert \nabla\overline
{u}\right\vert^{\ast}=\left\vert \nabla
{u}\right\vert^{\ast}$ and%

\begin{equation}
J\left(  u\right)  \geq\frac{%
{\displaystyle\int_{\Omega^{\#}}}
\left\vert \nabla\overline{u}\right\vert ^{2}dx-\dfrac{\left(  N-2\right)
^{2}}{4}%
{\displaystyle\int_{\Omega^{\#}}}
\dfrac{\overline{u}^{2}}{\left\vert x\right\vert ^{2}}dx}{\left\Vert
|\nabla\overline{u}|\right\Vert_{q}^{2}}.
\label {dec}
\end{equation}

\end{theorem}

{\normalsize \noindent\textbf{Proof .} In order to prove the theorem, it is enough to show that there exist a spherically symmetric decreasing function
$\overline{u}$ as in the statement of the theorem, such that $\|u\|_{L(2^{*},2)}\leq\|\overline{u}\|_{L(2^{*},2)}$. A result of this type, when the norm involved is the $L^{q}$ norm,
{\normalsize {with }$1\leq q\leq2N/\left(  N-2\right)  ,$ is
contained in \cite{ALT}. Other related results are also contained
in \cite{Cianchi} and \cite{Talenti}. }\\Set $f_0=\left\vert \nabla u\right\vert^{\ast}$.
By a result due to \cite{Giarrusso} (see also  \cite{ALT}) we get that
\begin{equation}
u^{\ast}\left(  s\right)  \leq\frac{1}{N\omega_{N}^{\frac{1}{N}}}\int_{s}%
^{|\Omega|}\frac{F\left(  t\right)  }{t^{1-\frac{1}{N}}}dt\text{, \
\thinspace }\forall s\in\left[  0,|\Omega|\right] \label{eq.1}  ,
\end{equation}
for some non negative function $F\in L^{2}\left(  0,|\Omega|\right)  $,
$F\prec f_{0}.$  The right-hand side of
(\ref{eq.1}) is the unique
spherically symmetric decreasing solution to the problem%
\[
\left\{
\begin{array}
[c]{l}%
\left\vert \nabla g\right\vert =F(\omega_{N}\left\vert
x\right\vert^{N})  \text{ \ in }\Omega^{\#}\\
\\
g=0\text{ \ on }\partial\Omega^{\#}.
\end{array}
\right.
\]
i.e.
\[
g\left(  \left\vert x\right\vert \right)  =\int_{\left\vert
x\right\vert }^{\left(  \frac{|\Omega|}{\omega_{N}}\right)
^{\frac{1}{N}}}F\left(  \omega _{N}s^{N}\right)  ds.
\]
The above result allows us to say that
\begin{equation}
\|u\|_{L\left(2^{*},2\right)\left(\Omega\right)}\leq \|g\|_{L\left(2^{*},2\right)(\Omega^{\#})}.\label{dis.1}
\end{equation}
At this point the proof is based on a duality argument, an approach that is different from the one shown in \cite{ALT}.
\noindent For each  $F\in L^{2}\left(  0,|\Omega|\right)  $, $F\prec f_{0},$ we define
\[
I\left(  F\right)  =\|g\|_{L\left(2^{*},2\right)\left(\Omega^{\#}\right)}=\left[  \int_{0}^{|\Omega|}\left(
\frac{1}{N\omega_{N}^{1/N}}\int_{s}^{|\Omega|}%
\frac{F\left(  t\right)  }{t^{1-1/N}}dt\right)
^{2}s^{-2/n}ds\right] ^{1/2}.
\]
On the other hand, since the dual of $L\left(  2^{\ast},2\right)
\left(\Omega^{\#}\right)  $ is the Lorentz space $L\left(  \frac{2N}{N+2},2\right)  \left(  \Omega^{\#}\right)
$, by definition of the norm in the dual space, for a fixed $F\in L^{2}\left(  0,|\Omega|\right)  $, $F\prec f_{0},$
we find%
\begin{equation}
I\left(  F\right)  =\max_{\substack{\phi\in L(\frac{2N}{N+2},2)\\\left\Vert \phi\right\Vert
_{\frac{2N}{N+2},2}=1}}\frac
{1}{N\omega_{N}^{1/N}}\int_{0}^{|\Omega|}\phi\left( s\right) \left(
\int_{s}^{|\Omega|}\frac{F\left(
t\right)  }{t^{1-1/N}}dt\right)  ds. \label{eq.11}%
\end{equation}
Let $\phi$ be a function for which the maximum in (\ref{eq.11}) is
attained.
Integrating by parts we obtain%
\[
\int_{0}^{|\Omega|}\phi\left(  s\right)  \left(  \int
_{s}^{|\Omega|}\frac{F\left(  t\right)  }{t^{1-1/N}%
}dt\right)  ds=\int
_{0}^{|\Omega|}\frac{F\left(  s\right)  }{s^{1-1/N}%
}\left[  \int_{0}^{s}\phi\left(  t\right)  dt\right]  ds.
\]
Hence%
\[
I\left(  F\right)  =\int_{0}^{|\Omega|}\frac{F\left( s\right)
}{N\omega_{N}^{1/N}s^{1-1/N}}\left[ \int_{0}^{s}\phi\left(
t\right)  dt\right]  ds.
\]
Let%
\[
\psi\left(  s\right)  =\frac{1}{N\omega_{N}^{1/N}s^{1-1/N}}\int_{0}^{s}%
\phi\left(  t\right)  dt,
\]
then $\psi\in L^{2}\left(  0, |\Omega|\right)  $. Indeed%
\begin{align*}
\left(  \int_{0}^{|\Omega|}\left[  \psi\left(  s\right) \right]
^{2}ds\right)  ^{1/2}  &  =\frac{1}{N\omega_{N}^{1/N}}\left(
\int_{0}^{|\Omega|} s^{-2+\frac{2}{N}}\left(  \int_{0}%
^{s}\phi\left(  t\right)  dt\right)  ^{2}ds\right)  ^{1/2}\\
&  =\frac{1}{N\omega_{N}^{1/N}}\left(
\int_{0}^{|\Omega|}s^{\frac{2}{N}}\left(
\frac{1}{s}\int_{0}^{s}\phi\left( t\right) dt\right) ^{2}ds\right)
^{1/2}\leq \frac{1}{N\omega_{N}^{1/N}}\left\Vert \phi\right\Vert
_{\frac{2N}{N+2},2}.
\end{align*}
Since $F\in L^{2}\left(0,|\Omega|\right)$ and  $F\prec f_{0}
$, we can find a sequence (see \cite{AT}, \cite{Mossino})
\[
\left\{  f_{n}\right\}  _{n\in\mathbb{N}}\subseteq L^{2}\left(
0,|\Omega|\right) \text{, }f_{n}^{\ast}=f_{0}
\] such that
$f_{n}\underset{n}{\rightharpoonup}F$ weak in $L^{2}\left(0,|\Omega|\right)$, then%
\[
I\left(  F\right)
=\lim_{n\rightarrow\infty}\int_{0}^{|\Omega|}f_{n}\left(  s\right)
\psi\left(  s\right) ds.
\]
\noindent By Hardy-Littlewood inequality
\[
\int_{0}^{|\Omega|}f_{n}\left(  s\right) \psi\left( s\right)
ds\leq\int_{0}^{|\Omega|}{f}_{0}\left(  s\right)  \psi^*\left(
s\right)ds,
\]
moreover (as it is shown in \cite {Chong}),  it is possible to construct a rearrangement $\overline{f}_{\psi}\in L^{2}\left(
0,|\Omega|\right)  $ of $f_{0}$ such that
\[
\int_{0}^{|\Omega|}f_{0}\left(  s\right) \psi^*\left( s\right)
ds=\int_{0}^{|\Omega|}\overline{f}_{\psi }\left(  s\right)  \psi\left(
s\right)ds  .
\]
The function  $\overline{f}_{\psi}$ is obtained by taking a sort of mean value of $f_0$, it essentially can be expressed in terms of the mean value operator introduced in  \cite{Mossino} and is connected with the notion of pseudo-rerrangement or relative rearrangement (see also  \cite{ALT}).

Therefore for any $F\in L^{2}\left(  0,|\Omega|\right)$ such that  $F\prec
f_{0}$ we have%
\begin{equation}
I\left(  F\right)  \leq\int_{0}^{|\Omega|}\overline{f}_{\psi}\left(
s\right) \psi\left(  s\right)  ds.\label{uno}
\end{equation}
Recalling the definition of $\psi$, an integration by parts allows
us to get that

\begin{equation}
\left.
\begin{array}
[c]{c}%
{\displaystyle\int_{0}^{|\Omega|}}
\overline{f}_{\psi}\left(  s\right)  \psi\left(  s\right)  ds=%
{\displaystyle\int_{0}^{|\Omega|}}
-\dfrac{d}{ds}\left(
{\displaystyle\int_{s}^{|\Omega|}}
\dfrac{\overline{f}_{\psi}\left(  t\right)  }{N\omega_{N}^{\frac{1}{N}%
}t^{1-\frac{1}{N}}}dt\right)  \left(
{\displaystyle\int_{0}^{s}}
\phi\left(  t\right)  dt\right)  ds\\
=\dfrac{1}{N\omega_{N}^{\frac{1}{N}}}%
{\displaystyle\int_{0}^{|\Omega|}}
\phi\left(  s\right)  \left(
{\displaystyle\int_{s}^{|\Omega|}}
\dfrac{\overline{f}_{\psi}\left(  t\right)  }{t^{1-\frac{1}{N}}}dt\right)
ds=I\left(  \overline{f}_{\psi}\right)
\end{array}
\right.\label{final}
\end{equation}

Hence setting
\[
\overline{u}\left(  \left\vert x\right\vert \right)
=\int_{\left\vert
x\right\vert }^{\left(  \frac{|\Omega|}{\omega_{N}}\right)  ^{\frac{1}{N}}}%
\overline{f}_{\psi}\left(  \omega_{N}s ^{N}\right)ds,
\]
by (\ref{dis.1}), (\ref{uno}), (\ref{final}) we have found a spherically symmetric decreasing function defined on $\Omega^{\#}$ such that $|\nabla \overline{u}|^*=f_0$ and for which

\[
\|u\|_{L\left(2^{*},2\right)\left(\Omega\right)}\leq \|\overline{u}\|_{L\left(2^{*},2\right)\left(\Omega^{\#}\right)},
\]
that is inequality (\ref {dec}) holds.$\qedsymbol$

{\normalsize \bigskip Now if we come back to the problem
concerning the calculation of the best constant $C$ in
(\ref{eq.20}), i.e.
the infimum%
\begin{equation}
C=\inf_{\substack{u\in H_{0}^{1}\left(  \Omega\right)  \\u\not \equiv 0}%
}\frac{\displaystyle%
{\displaystyle\int_{\Omega}}
\left\vert \nabla u\right\vert ^{2}-\dfrac{\left(  N-2\right)  ^{2}}{4}%
{\displaystyle\int\limits_{\Omega}}
\dfrac{u^{2}}{\left\vert x\right\vert ^{2}}dx}{\left\Vert\left \vert \nabla
u\right\vert\right\Vert_{q}^{2}}, \label{eq.12}%
\end{equation}
theorem 3 allows us to restrict our
attention to the
class of spherically symmetric decreasing functions $u\in H_{0}^{1}(\Omega^{\#}%
)$.

Indeed, if we reduce the study to spherically symmetric decreasing functions $u\in H_{0}^{1}(B_{R})$
with $R=\left(\left\vert\Omega
\right\vert/ \omega_{N}\right)^{1/N}$ , integrating by parts we have%
\[
\int_{B_{R}}\left\vert \nabla u\right\vert dx=-N\omega_{N}\int_{0}%
^{R}u^{\prime}\left(  r\right)  r^{N-1}dr=N\omega_{N}\left[  \lim
_{r\rightarrow0}u\left(  r\right)  r^{N-1}+\left(  N-1\right)  \int_{0}%
^{R}u\left(  r\right)  r^{N-2}dr\right]  .
\]
Since $H_{0}^{1}$ is imbedded $L\left(  2^{\ast},\infty\right)  $ we find that, for a
suitable constant $k$%
\[
u\left(  r\right)  \leq kr^{-\frac{N-2}{2}},
\]
from which it follows%
\[
\lim_{r\rightarrow0}u\left(  r\right)  r^{N-1}=0.
\]
Hence%
\[
\displaystyle\int_{B_{R}}\left\vert \nabla u\right\vert dx=\left(
N-1\right) \int_{B_{R}}\frac{u}{\left\vert x\right\vert }dx=\left(N-1\right)\omega_{N}^{\frac{1}{N}}\left\Vert u\right\Vert_{\frac{N}{N-1},1}.
\]
Applying theorem 2 with $p=N/(N-1)$ and hence $C(p,|\Omega|)=N^{2}\omega_{N}^{2/N}%
/(4|\Omega|)$, we deduce%

\begin{align*}
C  & =\frac{1}{\omega_{N}^{2/N}\left(  N-1\right)
^{2}}\inf_{\substack{u\in
H_{0}^{1}(\Omega^{\#})\\u=u^{\#}\\u\not \equiv 0}}\dfrac{%
{\displaystyle\int_{\Omega^{\#}}}
\left\vert \nabla u\right\vert ^{2}dx-\dfrac{\left(  N-2\right)  }{4}%
{\displaystyle\int_{\Omega^{\#}}}
\dfrac{u^{2}}{\left\vert x\right\vert ^{2}}dx}{\left\Vert
u\right\Vert
_{\frac{N}{N-1},1}^{2}}\\
& =\frac{1}{4\left\vert \Omega\right\vert }\left(
\frac{N}{N-1}\right)  ^{2}.
\end{align*}
In conclusion we have
\begin{theorem}
{\normalsize  Let $\Omega$ be a bounded open subset of
$\mathbb{R}$$^{N}$ containing the origin, $N>2$. If $q=1$, the optimal value of the
constant in inequality (\ref{eq.20}) is given by
\[
C(|\Omega|)=\frac{1}{4\omega_N\left\vert \Omega\right\vert }\left(  \frac{N}%
{N-1}\right)  ^{2}.
\]
}
\end{theorem}

\bigskip
\bigskip
\bigskip

It is clear that it is always possible to reduce the calculation of the best constant $C$ to a spherically symmetric situation in inequalities of the type (5) that involve the Lorents norms $\left\Vert\left \vert \nabla
u\right\vert\right\Vert_{p,q}$ with $1\leq q<p<2$ or $\left\Vert\left \vert \nabla
u\right\vert\right\Vert_{p,1}$ with $0<p<1$. In particular, the method used in the proof of theorem 4 allows us to
treat inequalities of the type (\ref{eq.20}) when the right hand side is
\[
\int_{\Omega^{\#}}\left\vert \nabla u\right\vert ^{\#}\left\vert
x\right\vert ^{\alpha}dx\text{, \ with }\alpha>0,
\]
that is essentially a norm of $|\nabla u|$ in a Lorentz
space $L(p,1)$ with $0<p=\frac {N}{N+\alpha}<1$. Indeed, also in this case, by an integration by parts, the problem is reduced to the study of an inequality of the type (\ref{caso31}). In fact, using the
Hardy-Littlewood inequality, since the function $|x|^\alpha$ is
increasing we find
\begin{equation}
\int_{\Omega^{\#}}\left\vert \nabla u\right\vert ^{\#}\left\vert
x\right\vert ^{\alpha}dx   \leq\int_{\Omega^{\#}}\left\vert
\nabla u\right\vert \left\vert
x\right\vert ^{\alpha}dx=\left(  N+\alpha-1\right)  \int_{\Omega^{\#}}u\left\vert
x\right\vert ^{\alpha-1}dx\label{rem1}\\
\end{equation}
for all $u\in H_{0}^{1}\left(\Omega^{\#}\right)$ such that
$u=u^{\#}$.
More precisely, the following theorem can be proved:

\begin{theorem} Let $\Omega$ be a bounded open subset of
$\mathbb{R}$$^{N}$ containing the origin, $N>2$. If $q=1$ and $\max\left\{\frac{2N}{3N-2},\frac{N}{N+1}\right\}<p<1$, the optimal value of the
constant in  the inequality

\begin{equation}
\int_{\Omega}\left\vert \nabla u\right\vert ^{2}dx-\frac{\left(
N-2\right) }{4}^{2}\int_{\Omega}\frac{u^{2}}{\left\vert
x\right\vert ^{2}}dx\geq
C\left\Vert \left\vert\nabla u\right\vert\right \Vert _{p,1} ^{2}\text{, \ \ }\forall u\in H_{0}%
^{1}\left(  \Omega\right) , %
\end{equation}
\noindent is given by

\[
C(p,|\Omega|)=\frac{\left(\frac{2-p}{p}\right)^3}{4\left\vert \Omega\right\vert ^{2/p-1}}\left(  \frac{Np}%
{N-p}\right)  ^{2}.
\]

\end{theorem}

\end{document}